\documentclass[11pt]{article}
\usepackage{amssymb}
\usepackage{graphicx}
\usepackage{amsmath}
\usepackage{amssymb}
\usepackage{graphicx}
\usepackage{amsmath}

\newcommand{\Z}{\mathbb {Z}}
\newcommand{\Q}{\mathbb {Q}}

\newcommand{\R}{\mathbb {R}}
\newcommand{ \Spec}{\mbox{\textup{Spec }}}

\newcommand{ \Pic}{\mbox{\textup{Pic}}}
\newcommand{ \Rad}{\mbox{\textup{Rad}}}

\newcommand{ \Rees}{\mbox{\textup{Rees }}}

\newcommand{ \Char}{\mbox{\textup{char }}}

\def\inbar{\,\vrule height1.5ex width.4pt depth0pt}
\def\inbar{\,\vrule height1.5ex width.4pt depth0pt}
\def\IB{\relax{\rm I\kern-.18em B}}
\def\IC{\relax\hbox{$\inbar\kern-.3em{\rm C}$}}
\def\ID{\relax{\rm I\kern-.18em D}}
\def\IE{\relax{\rm I\kern-.18em E}}
\def\IF{\relax{\rm I\kern-.18em F}}
\def\IG{\relax\hbox{$\inbar\kern-.3em{\rm G}$}}
\def\IH{\relax{\rm I\kern-.18em H}}
\def\II{\relax{\rm I\kern-.18em I}}
\def\IK{\relax{\rm I\kern-.18em K}}
\def\IL{\relax{\rm I\kern-.18em L}}
\def\IM{\relax{\rm I\kern-.18em M}}
\def\IN{\relax{\rm I\kern-.18em N}}
\def\IO{\relax\hbox{$\inbar\kern-.3em{\rm O}$}}
\def\IP{\relax{\rm I\kern-.18em P}}
\def\IQ{\relax\hbox{$\inbar\kern-.3em{\rm Q}$}}
\def\IR{\relax{\rm I\kern-.18em R}}
\font\cmss=cmss10 \font\cmsss=cmss10 at 7pt
\def\IZ{\relax\ifmmode\mathchoice
{\hbox{\cmss Z\kern-.4em Z}}{\hbox{\cmss Z\kern-.4em Z}}
{\lower.9pt\hbox{\cmsss Z\kern-.4em Z}}
{\lower1.2pt\hbox{\cmsss
Z\kern-.4em Z}}\else{\cmss Z\kern-.4em Z}\fi}
\def\IGa{\relax\hbox{${\rm I}\kern-.18em\Gamma$}}
\def\IPi{\relax\hbox{${\rm I}\kern-.18em\Pi$}}
\def\ITh{\relax\hbox{$\inbar\kern-.3em\Theta$}}
\def\IOm{\relax\hbox{$\inbar\kern-3.00pt\Omega$}}

\begin{document}

\baselineskip 20pt
\pagenumbering{arabic}
\pagestyle{plain}

\newtheorem{defi}{Definition}[section]
\newtheorem{theo}[defi]{Theorem}
\newtheorem{lemm}[defi]{Lemma}
\newtheorem{prop}[defi]{Proposition}
\newtheorem{note}[defi]{Note}
\newtheorem{nota}[defi]{Notation}
\newtheorem{exam}[defi]{Example}
\newtheorem{coro}[defi]{Corollary}
\newtheorem{rema}[defi]{Remark}
\newtheorem{cons}[defi]{Construction}
\newtheorem{ques}[defi]{Question}
\newtheorem{questions}[defi]{Questions}

\newcommand{\abs}{${\bar A}^*}
\newcommand{\qbs}{${\bar Q}^*}
\newcommand{\be}{\begin{enumerate}}
\newcommand{\ee}{\end{enumerate}}
\newcommand{\fany}{\rm for\ \ any\ \ }

\newcommand{\rb}{\overline{R}}
\newcommand{\mt}{\overline{M}}
\newcommand{\nt}{\overline{N}}
\newcommand{\nb}{\widetilde{N}}
\newcommand{\mb}{\widetilde{M}}
\newcommand{\m}{\bf {m}}


\def\cm{Cohen-Macaulay}
\def\wrt{with respect to\ }
\def\pni{\par\noindent}
\def\wma{we may assume without loss of generality that\ }
\def\Wma{We may assume without loss of generality that\ }
\def\ets{it suffices to show that\ }
\def\bwoc{by way of contradiction}
\def\iff{if and only if\ }
\def\st{such that\ }
\def\fg{finitely generated}

\def\R{\Cal R}
\def\G{\Cal G}
\def\F{\Cal F}
\def\a{\goth a}

\def\M{\Cal M}
\def\N{\goth N}
\def\P{\Cal P}
\def\I{\Cal I}
\def\J{\Cal J}
\def\Q{\Cal Q}
\def\p{\mathbf p}

\def\isom{\thinspace \cong\thinspace}
\def\rtar{\rightarrow}
\def\rta{\rightarrow}
\def\l{\lambda}
\def\d{\Delta}

\def\alert#1{\smallskip{\hskip\parindent\vrule%
\vbox{\advance\hsize-2\parindent\hrule\smallskip\parindent.4\parindent%
\narrower\noindent#1\smallskip\hrule}\vrule\hfill}\smallskip}

\title{
PROJECTIVELY FULL IDEALS IN \\ NOETHERIAN RINGS, A SURVEY}
\author{
Catalin Ciuperca,
William Heinzer,
Jack Ratliff and David Rush}

\date{\today}
\maketitle

\begin{abstract}
We discuss projective equivalence of ideals in  Noetherian rings
and the existence or failure of existence of projectively full ideals. 
We describe 
connections with the Rees valuations and  Rees integers of an ideal,
and consider the question of whether improvements can be made by
passing to an integral extension ring.

\end{abstract}

\section{Definitions, summary and examples.}

Let $R$ be a Noetherian ring and let $I \neq R$
be a regular ideal of $R$. (So $I$ contains an element
with zero annihilator.)

\begin{defi} {\em
An ideal $J$ of $R$ is {\bf projectively equivalent} to $I$ if there
exist positive integers $m$ and $n$ such that $I^m$ and $J^n$ have
the same integral closure, i.e., $(I^m)_a = (J^n)_a$.}
\end{defi}

Projective equivalence is an equivalence relation on the regular
proper ideals of $R$.  Let $\mathbf P(I)$ denote the set of
integrally closed ideals projectively equivalent to $I$.

\begin{rema} {\em
The set $\mathbf P(I)$ is discrete and linearly ordered with respect
to inclusion. Moreover,  $J$ and $K$ in $\mathbf P(I)$ and $m$ and
$n$ positive integers implies $(J^mK^n)_a \in \mathbf P(I)$. Thus
there is naturally associated to $I$ and $\mathbf P(I)$ a unique
subsemigroup $S(I)$ of the additive semigroup of nonnegative
integers $\IN_0$  such that $S(I)$ contains all sufficiently large
integers. A semigroup having these properties is called a {\bf
numerical semigroup}. It is observed in \cite[Remark 3.11]{CHRR2}
that every numerical semigroup $S$ is realizable as $S(M)$, where
$(R,M)$ is a local domain of altitude one.

}
\end{rema}

\begin{defi} {\em
The set $\mathbf P(I)$ is said to be {\bf
projectively full}
if $S(I) = \IN_0$, or equivalently,  if every element of $\mathbf
P(I)$ is the integral
closure of a power of the largest element $K$ of $\mathbf P(I)$,
i.e., every
element of $\mathbf P(I)$ has the form $(K^n)_a$, for some
positive integer $n$.
If this holds, then each ideal $J$ in $R$ such that $J_a = K$
is said to be {\bf projectively
full}.

}
\end{defi}

\medskip

\noindent
Our main goals may be listed as follows:

\noindent
{\bf 1}.  Describe ideals $K$ that are projectively full.

\medskip

\noindent
{\bf 2}.  Describe ideals $I$ such that $\mathbf P(I)$ is
projectively full.

\medskip

\noindent
{\bf 3}.  Determine if improvements can be obtained by passing to an
integral extension ring.

\medskip

In connection with the third main goal, the following theorem is
proved in \cite{CHRR3}.

\begin{theo}
If $R$  contains the field of rational numbers, then
 there exists a finite free integral extension ring $A$  of
$R$ such that
$\mathbf P(IA)$  is projectively full; and if $R$ is an integral
domain, then
there also exists a finite integral extension domain $B$  of $R$
such
that $\mathbf P(IB)$  is projectively full.

\end{theo}

A strengthened version of Theorem 1.4 is given in
Corollary \ref{great} below.

Examples \ref{1.5} and \ref{1.6} illustrate aspects of the
projectively full property.

\begin{exam} \label{1.5}
{\em
Let $(R,M)$ be a Noetherian local ring having the property
$$a \in M^i \setminus M^{i+1}  \text{ and } b \in M^j
\setminus M^{j+1}  \implies ab \not\in M^{i+j+1},
$$
then $M$ is projectively full. Thus if the associated graded ring
$$
\mathbf G(R, M) = R/M \oplus M/M^2 \oplus \cdots \oplus
M^n/M^{n+1} \oplus \cdots
$$
is an integral domain, then $M$ is projectively full.
}
\end{exam}

\begin{exam} \label{1.6}
{\em Let $x$ be an indeterminate over the field $F$  and let $R =
F[[x^2, x^3]]$. Then the maximal ideal $M = (x^2, x^3)R$ is not
projectively full. The numerical semigroup $S(M)$ is generated by
$2$ and $3$. }
\end{exam}

Notice that with $R = F[[x^2, x^3]]$ as in Example \ref{1.6}, the
ring  $R$ is not integrally closed. Things improve by passing to the
integral extension ring $R[x] = F[[x]]$. Indeed, for each regular
proper ideal $I$ of $R$, $\mathbf P(IR[x])$ is projectively full.

\medskip

Example \ref{1.7} is established in \cite{CHRR2} to demonstrate
existence of a normal local domain $(R,M)$ such that $M$ is not
projectively full.

\begin{exam} \label{1.7}
{\em
Let $F$ be a field in which $2$ and $3$ are units and let $x,
y, z, w$ be variables. Let
$$ R_0 = F[x, y]_{(x, y)}\quad \text{  and  } \quad  R = \frac{R_0[z, w]}{(z^2 - x^3 - y^3, \,\, w^2
- x^3 + y^3)}.
$$
Then $R$ is a normal local domain of altitude\footnote{We are
following Nagata \cite{N2} and  using altitude for what is often now
termed dimension or Krull dimension.} two with maximal ideal $M =
(x,y, z, w)R$, and $M$ is not projectively full.

}
\end{exam}

\section{Some history.}
The concept of projective equivalence of ideals and the study of
ideals projectively equivalent to $I$ was introduced by Samuel in
\cite{Samuel} and further developed by Nagata in \cite{Nagata}.
Rather than `projectively equivalent', Hartmut  G\"ohner uses the
term `asymptotically equivalent' in \cite{Gohner}. G\"ohner mentions
that the expression `projective asymptotic equivalence' is used by
David Rees in \cite{Rees.val.assoc.to.ideals.II} and by H. T. Muhly
in \cite{Muhly}.

Let $I$ be a regular proper ideal of a Noetherian ring $R$. For each
$x \in R$,   let $v_I(x)$ $=$ $\max \{k \in \IN  \mid x \in I^k\}$,
and define
  $$\overline v _I (x)=\lim _{k \rightarrow \infty}
(\frac{v_I(x^k)}{k}).$$
Rees established that:

\noindent
{\bf (a) }  $\overline v _I(x)$  is well defined.

\noindent
{\bf (b) }  For each  $k$ $\in$ $\IN$ and  $x$ $\in$ $R$, $\overline v _I
(x)$ $\ge$ $k $   if and only if  $x \in (I^k)_a$.

\noindent
{\bf (c) } There exist discrete valuations  $v_1,\dots,v_n$  defined on
$R$,  (with values in $\IN \cup \infty$ ) and
positive integers $e_1,\dots,e_n$ such that, for each $x$ $\in$
$R$,
$\overline v _I(x) = \min \{\frac{v_i(x)}{e_i} \mid i =
1,\dots,n\}$.

In (c), we say
that $v$ is a valuation on $R$ if $\{x \in R \, | \, v(x) = \infty \}$ is
a prime ideal $P$, $v(x) = v(y)$ if $x + P = y + P$, and the function $v$
induces on $R/P$ is a valuation.

For simplicity, assume $R$ is a Noetherian integral domain with
field
of fractions $F$. Let  $t$  be an indeterminate.  The {\bf Rees
ring } of $R$
with respect to $I$ is the graded subring
$$
\mathbf R = R[t^{-1}, It] = \bigoplus_{n \in \Z}I^nt^n
$$
of the Laurent polynomial ring $R[t^{-1}, t]$. The integral
closure
$B$ of $R[t^{-1}, It]$ is a Krull domain, and $B_P$ is a DVR for
each
minimal prime $P$ of $t^{-1}B$. The set $\Rees I$ of {\bf Rees
valuation
rings} of $I$ is  precisely the set of rings  $V = B_P \cap F$,
where $P$ is
a minimal prime of $t^{-1}B$.

\begin{defi} {\em
Let $(V_1, N_1), \ldots, (V_n, N_n)$ be the Rees valuation rings of
$I$. The integers $(e_1, \ldots, e_n)$, where $IV_i = N_i^{e_i}$,
are the {\bf Rees integers } of $I$.  (These are the same
$e_i$  as in (c) above.)

}
\end{defi}

A sufficient condition for $I$ to be projectively full is that
$\gcd(e_1, \ldots, e_n) = 1$. However, this is not a necessary
condition as we demonstrate in Example \ref{3.1}.

\section{Projective equivalence and Rees valuation rings.}

\begin{exam} \label{3.1}
{\em Let $(R,M)$ be a regular local ring of altitude two  with $M = (x,
y)R$, and let $e > 1$ be an integer. The ideal $I = (x, y^e)R$ is
integrally closed with unique Rees valuation ring $V =
R[x/y^e]_{MR[x/y^e]}$. The integrally closed ideals projectively
equivalent to $I$ are precisely the powers $I^n$ of $I$. Thus $I$ is
projectively full. The maximal ideal of $V$ is $N = yV$ and $IV =
N^e$, so the gcd  of the Rees integers of $I$ is $e > 1$. }
\end{exam}

If $I$ and $J$ are projectively equivalent, then $\Rees I = \Rees
J$. The converse holds if $I$ or $J$ has only one Rees valuation
ring, but fails in general. Steve McAdam, Jack Ratliff and Judy
Sally show in \cite{MRS} that if $I$ and $J$ are projectively
equivalent, then the Rees integers of $I$ and $J$ are proportional.
It is observed in \cite{CHRR} that the converse also holds:  if
$\Rees I = \Rees J$ and the Rees integers of $I$ and $J$ are
proportional, then $I$ and $J$ are projectively equivalent.

\begin{exam} \label{3.2}
{\em Assume $(R,M)$ is a regular local ring of altitude two.
Zariski's theory of unique factorization of
complete (= integrally closed) ideals of $R$ as  finite products of
simple complete ideals implies $\mathbf P(I)$ is projectively full
for each nonzero proper ideal $I$ of $R$. The ideal $I$ has a unique
Rees valuation ring if and only if $I$ is a power of a simple
complete ideal. If $I$ factors as
$I = I_1^{f_1} \cdots I_n^{f_n}$,
where $I_1, \ldots, I_n$ are distinct simple complete ideals, then
$I$ is
projectively full if and only if $\gcd(f_1, \ldots, f_n) = 1$.

How do the integers $f_1, \ldots, f_n$ relate to the Rees integers
of $I$?
The simple complete $M$-primary ideals of $R$ are in one-to-one
correspondence with the prime divisors birationally dominating
$R$.
Thus the Rees valuation rings of $I$ are  $(V_1, N_1), \ldots,
(V_n, N_n)$, where
$(V_j, N_j)$ corresponds to $I_j$. If $I_jV_j = N_j^{c_j}$, then
the Rees integers
of $I$ are \,\,  $e_1 = c_1f_1, \, \, \ldots,\,\,  e_n = c_nf_n$.
}
\end{exam}

Example \ref{3.3} illustrates these concepts in a specific
situation.

\begin{exam} \label{3.3}
{\em Assume $(R,M)$ is a regular local ring of altitude two with $M =
(x, y)R$ as in Example \ref{3.1}. Let $I = (x^2, xy^2, y^3)R$.
Notice that $J = (x^2, y^3)R$ is a reduction of $I$ and $JI = I^2$,
so the reduction number $r_J(I) = 1$. Let
$$V = R[\frac{xy^2}{x^2}, \frac{y^3}{x^2}]_{MR[\frac{xy^2}{x^2},
\frac{y^3}{x^2}]} = R[\frac{y^2}{x},
\frac{y^3}{x^2}]_{MR[\frac{y^2}{x}, \frac{y^3}{x^2}]}.
$$
One sees that $V$ is a valuation ring with maximal ideal
$N = (y^2/x)V$, and $I$ is a simple complete ideal. The ideals of
$R$ that are contracted from $V$ descend as follows:
$M = N^2 \cap R \supsetneq (x, y^2)R = N^3 \cap R \supsetneq M^2 =
N^4 \cap R \supsetneq (x, y^2)M$
$= N^5 \cap R \supsetneq I = N^6 \cap R $. The ideals in $\mathbf
P(I)$ are precisely the
ideals $I^m = N^{6m} \cap R$, for $m \in \IN$.
}
\end{exam}

\section{The Rees integers and the associated graded ring.}

The following result is established in \cite[Example 5.1]{CHRR3}.

\begin{prop} \label{4.1}
If $\mathbf G(R,I) = R/I \oplus I/I^2 \oplus \cdots \oplus
I^n/I^{n+1} \oplus \cdots$ has a minimal prime $z$ such that $z$ is
its own $z$-primary component of $(0)$, then $I$ has a Rees integer
equal to one. Therefore $I$ is projectively full. In particular, if
$\mathbf G(R,I)$ is reduced, then $I$ is projectively full.

\end{prop}

More can be said about the Rees integers of $I$ by using the Rees
ring $\mathbf R = R[t^{-1}, It]$  and the identification $\mathbf
G(R,I) = \mathbf R/t^{-1}\mathbf R$. \,\, Let $\mathbf R'$ denote
the integral closure of $\mathbf R$. Proposition \ref{4.2} is  from
\cite[Proposition 4.1]{CHRR3}.

\begin{prop} \label{4.2}
The ideal $I$ has a Rees integer equal to one if and only if
$t^{-1}\mathbf R$ has a minimal prime $p$ such that $t^{-1}\mathbf
R'_p = p\mathbf R'_p$.

\end{prop}

Let $(R,M)$ be a normal local domain and let $\mathbf G(R,M)$ denote
the associated graded ring of $R$ with respect to $M$. Proposition
\ref{4.1} implies that $M$ is projectively full if $\mathbf G(R,M)$
is reduced. Example \ref{4.3} shows that the converse of this
statement is not true.

\begin{exam} \label{4.3}
{\em
An example of a normal local domain $(R,M)$ of altitude two such
that $M$ is projectively full and the associated graded ring
$\mathbf G(R,M)$ is not reduced. Let $F$ be an algebraically closed
field with $\Char F = 0$,
 and let $R_0$ be a regular local domain of altitude two with
 maximal
 ideal $M_0 = (x, y)R_0$ and coefficient field $F$, e.g.,
 $R_0 = F[x, y]_{(x, y)}$, or $R_0 = F[[x, y]]$.
 Then
$V_0 = R_0[y/x]_{xR_0[y/x]}$
is the unique Rees valuation ring of $M_0$. Let
 $R = R_0[z]$,  where   $z^2 =  x^3 + y^j$
with   $j \ge 3$.
 It is readily checked that
 $R$ is a normal local domain  of altitude two with maximal
 ideal $M =
 (x, y, z)R$. Moreover, the image of $z$ in $M/M^2$ is a
nonzero nilpotent element of the associated graded ring  $\mathbf G(R,M)$,
so $\mathbf G(R,M)$ is not reduced.

Notice
 that $I = (x, y)R$ is a reduction of $M$
 since $z$ is integral over $I$. It follows that
every Rees valuation ring of $M$ is an   extension  of
$V_0$.
Let $V$ be a Rees valuation ring of $M$ and let $v$ denote the
normalized valuation with value group $\Z$ corresponding to $V$.
Then $v(x) = v(y)$ and the image of $y/x$ in the residue field of
$V$
is transcendental over $F$. Since $ z^2 = x^3 + y^j$ \,\,and
\,\,$j \ge 3$,
 we have
$$
2v(z) = v(z^2) = v(x^3 + y^j) = 3v(x). $$ We conclude that  $v(x) =
2$ and $v(z) = 3$. Therefore $V$ is ramified over $V_0$.  This
implies  that $V$ is the unique extension of $V_0$ in the quadratic
field extension generated by $z$ over the field of fractions of
$R_0$, and thus the unique Rees valuation ring of $M$.

For each positive integer $n$, let
$I_n = \{ r \in R \, | \, v(r) \ge n \}$. Thus $I_2 = M$. Since
$V$ is the unique
Rees valuation ring of $M$, we have $I_{2n} = (M^n)_a$ for each $n
\in \IN$.

To show $M$ is projectively full,  we prove that  $V$ is not the
unique Rees valuation ring of $I_{2n+1}$ for each $n \in \IN$.
Consider the inclusions
$$ M^2 \subseteq I_4 \subset (z, x^2, xy, y^2)R := J \subseteq I_3 \subset
I_2 = M. $$
Since $\lambda(M/M^2) = 3$ and since the images of $x$ and $y$ in
$M/M^2$
are $F$-linearly independent, \,\,  $J = I_3$ and $M^2 = I_4 =
(M^2)_a$.

Since $x^3 = z^2 - y^j$ and $j \ge 3$,
$L = (z, y^2)R$ is  a reduction of $I_3 = (z, x^2, xy, y^2)R$.
\,\, Indeed,
$(x^2)^3 \in L^3$ and $(xy)^3 \in L^3$ implies $x^2$ and $xy$ are
integral over $L$.
It follows that $V$ is not a Rees valuation of $I_3$,  for $zV
\ne y^2V$.

Consider $M^3 \subset I_3M \subseteq I_5 \subset I_4 = M^2$. Since
the images of $x^2, xy, y^2, xz, yz$ in $M^2/M^3$ are an $F$-basis,
it follows that $I_3M = I_5$ and $M^3 = (M^3)_a = I_6$. Proceeding
by induction, we assume $M^{n+1} = (M^{n+1})_a = I_{2n+2}$, and that
$I_{2n+1}$ is not projectively equivalent to $M$. Consider
$$  M^{n+2} \subset I_3M^n \subseteq I_{2n+3} \subset
I_{2n+2} = M^{n+1}.
$$

Since the images in
$M^{n+1}/M^{n+2}$
of $\{x^ay^b \, | \, a + b = n+1 \} \cup \{zx^ay^b \, | \, a + b =
n \}$ is an $F$-basis,
 $\lambda(M^{n+1}/M^{n+2}) = 2n+3$, and the inequalities
$\lambda(M^{n+1}/I_{2n+3}) \ge n+2$  and
$\lambda(I_3M^n/M^{n+2}) \ge n+1$
imply $I_3M^n = I_{2n+3}$  and  $M^{2n+2} =
(M^{2n+2})_a$.

Therefore the ideal $I_{2n+3}$ has a Rees valuation ring different
from $V$, and thus is not projectively equivalent to $M$. We
conclude that $M$ is projectively full. We have also shown
that $M$ is a normal ideal.
}
\end{exam}

\section{Questions and examples in altitude two.}

In Questions \ref{5.1} we list several questions that  interest  us
and that we hope may stimulate further research work in this area.

\begin{questions} \label{5.1} {\em
Let $(R,M)$ be a complete normal local domain of
altitude two.

\begin{enumerate}
\item  What are necessary and sufficient conditions in order that
$M$ be projectively full?

\item Let $I$ be an $M$-primary ideal. What are
necessary and sufficient conditions in order that $I$ be
projectively full?

\item What are necessary and sufficient conditions in order that
$\mathbf P(I)$  be projectively full for each $M$-primary  ideal $I$
of $R$?

\item What are necessary and sufficient conditions in order that
$\mathbf P(I)$  be  projectively full for each nonzero proper ideal
$I$ of $R$?

\end{enumerate}

}
\end{questions}

In Example \ref{5.2}, we present in more detail a family of examples
from \cite[Example 5.1]{CHRR2} that includes Example \ref{1.7}.

\begin{exam} \label{5.2}
{\em An example of a (complete)  normal local domain $(R,M)$ of
altitude two such that $M$ is not projectively full.  Let $F$ be an
algebraically closed  field,
 and let $R_0$ be a  regular local domain of altitude two  with  maximal
 ideal $M_0 = (x, y)R_0$ and
 coefficient field $F$, e.g.,
 $R_0 = F[x, y]_{(x, y)}$, or $R_0 = F[[x, y]]$.
 Let $k < i$ be relatively prime positive integers $\ge 2$,
 and let
$R = R_0[z, w]$,  where \, $z^k = x^i + y^i$ \,and \,  $w^k = x^i -
y^i$. Assume that $2$, $i$  and $k$ are units of $F$.
 It is readily checked that
 $R$ is a normal local domain of altitude two  with maximal
 ideal $M =
 (x, y, z,w)R$. Also $R$ is a free $R_0$-module of rank $k^2$.

With $R = R_0[z, w]$ as above, we
show  that $M = (x, y, z, w)R$ has a unique Rees
valuation ring
and that $M$ is not projectively full. Notice that  $L = (x, y)R$ is a
reduction of
$M$,  for  $z^k \in (x^i, y^i)R \subseteq L^k$   and
$w^k \in (x^i, y^i)R \subseteq L^k$ \, imply
$z$ and $w$ are integral over $L$. Thus each Rees valuation ring
$V$ of $M$
is an extension of the order valuation ring $V_0 =
R_0[y/x]_{xR_0[y/x]}$ of
 $R_0$. To show there exists a unique Rees valuation ring of $M$,
 we
 observe that if  $V$ is a Rees valuation ring of $M$, then
$V$ as an extension of $V_0$ ramifies of  degree
 $k$ and
 undergoes a residue field extension of degree at least  $k$.

 To show $V$ ramifies of  degree  $k$ over $V_0$, observe that
$kv(z) = v(z^k) = v(x^i + y^i) = iv(x) = iv(y)$ implies $v(z) = i$
and $v(x) = v(y) = k$. Similarly, $v(w) = i$. Let $N$ denote the
maximal ideal of $V$. We have $(z, w)V = N^i$ and $MV = (x, y)V =
N^k$. The residue field of $V_0$ is $F(\overline\tau)$, where
$\overline \tau$ is the image of $\tau = y/x$ and is transcendental
over $F$. Now $w/z$ is a unit of $V$ and
$$
(\frac{w}{z})^k  =
\frac{x^i - y^i}{x^i + y^i} = \frac{1 - \tau^i}{1 + \tau^i}.
$$  It
follows that the residue class of $w/z$ in $V/N$ is algebraic of
degree $k$ over $F(\overline \tau)$. \, This proves that $V$ is the
unique Rees valuation ring of $M$.

To show $M$ is not projectively
full, notice that $(z^k, w^k)R = (x^i + y^i, x^i - y^i)R$, and
since
$\Char F \ne  2$, \,\, $(x^i + y^i, \,  x^i - y^i)R = (x^i,
y^i)R$. \,\, Since $(x^i, y^i)R$
 is a reduction of $M^i$, we have $((z^k, w^k)R)_a =
(M^i)_a$. \, Also
$((z^k, w^k)R)_a = ((z, w)^kR)_a$. Therefore $(z, w)R$ and $M$ are
 projectively equivalent, so $V$ is the unique Rees valuation ring
of $(z, w)R$, and $((z, w)R)_a = N^i \cap R$. We have $M = N^k
\cap R$,
$M^nV = N^{nk}$ and $(M^n)_a = N^{nk} \cap R$, for each positive
integer $n$.
Since $i$ is not a multiple of $k$, $M$ is  not projectively full.
}
\end{exam}

\begin{rema} {\em
Let $(R,M)$ be a  normal local domain of altitude two.
 Muhly and Sakuma in
\cite{Muhly.Sakuma} define {\bf the filtration} $\{ A_i \mid i \in
\IN_0 \}$ {\bf associated to a prime $R$-divisor} $v$ by $A_0$ = $R$
and $A_i$ = $\{ x \in R \mid v(x) > v(A_{i-1} \}$ for $i \geq 1$.
The prime $R$-divisor $v$ is said to be {\bf Noetherian} if the
associated filtration $\{ A_i \mid i \geq 0 \}$ possesses a
Noetherian subfiltration, that is a subfiltration $\{ A_{{i_j}} \mid
j \geq 0\}$ such that the graded ring $\oplus_{j \geq 0}
A_{{i_j}}/A_{i_{j+1}}$ is Noetherian.
 Muhly and Sakuma say that a normal local domain $R$
{\bf satisfies condition N} if each prime divisor of the second kind
relative to $R$ is Noetherian \cite[First sentence of Section
4]{Muhly.Sakuma}. They show that if $R$ is analytically normal of
altitude two, then each prime divisor of the second kind relative to
$R$ is the unique valuation ring of some ideal $I$ of $R$ if and
only if $R$ satisfies condition N \cite[Proposition 2.2 and the
paragraph following it]{Muhly.Sakuma}. (It is shown in \cite{BPRR}
that some caution should be used in reading the given proof of their
needed result \cite[Proposition 2.1]{Muhly.Sakuma}.) Muhly and
Sakuma  also show in \cite{Muhly.Sakuma} that if an analytically
normal local domain $R$ of altitude two satisfies condition (N),
then there is
 a factorization theory for complete
$M$-primary ideals that extends some aspects of Zariski's
factorization theory for complete ideals in the case where $R$ is
regular. Muhly carries this work further in \cite{Muhly} by showing
that the  hypothesis that  $R$ is analytically normal  can be
weakened to $R$ is normal, and also that contrary to what was
conjectured in \cite{Muhly.Sakuma},  $R$ being analytically normal
does not imply $R$ satisfies  condition (N).

Assume $(R,M)$ is a complete normal local domain of altitude two.
G\"ohner in \cite{Gohner} proves that $R$ satisfies condition (N) if
$R$ is semifactorial, i.e., if $\Pic R$ is a torsion group, and
conjectures that the converse also holds. Cutkosky proves this
conjecture  in \cite{Cut}. This establishes an interesting
connection between the divisorial ideals of $R$ and the complete
$M$-primary ideals.

G\"ohner  \cite[Corollary 4.5]{Gohner} observes that if $(R,M)$ is a
complete normal local domain of altitude two and if $R/M$ is the
algebraic closure of a finite field, then $R$ is semifactorial and
therefore satisfies condition (N). By taking $F$ in Example
\ref{5.2} to be the algebraic closure of a finite field in which
$2$, $i$ and $k$ are units, we see that there exists a complete
normal local domain $(R,M)$ that satisfies condition (N) and has the
property that its maximal ideal $M$ is not projectively full.

}
\end{rema}

\begin{rema} {\em
Let $(R,M)$ be a  normal local domain of altitude two and let $I$ be
an  $M$-primary ideal. Recall that the multiplicity
$e(I) = e(I,R)$ of $I$  is a positive integer \cite[Section
14]{Mat}. In connection with Questions \ref{5.1}, Joe Lipman pointed
out to us that if $e(I)$ is squarefree, then $I$ is projectively
full.  For an ideal and its integral closure have the same
multiplicity, and for each positive integer $n$, $e(I^n) = n^2e(I)$
\cite[Formula 14.3]{Mat}. If $m$ and $n$ are  relatively prime
positive integers and $J$ is an ideal such that $(I^m)_a = (J^n)_a$,
then $m^2e(I) = n^2e(J)$, and the assumption that  $e(I)$ is
squarefree implies that $n = 1$. Therefore the hypothesis that
$e(I)$ is squarefree implies that $I$ is projectively full.

More generally, let $(R,M)$ be a local domain of altitude $d \ge 2$
and let $I$ be an $M$-primary ideal. For each positive integer $n$,
we have $e(I^n) = n^de(I)$. Assume that $m$ and $n$ are relatively
prime positive integers and that $J$ is an ideal such that $(I^m)_a
= (J^n)_a$. Then $m^de(I) = n^de(J)$. If $e(I)$ is free of $d$-th
powers, then $n = 1$. Therefore the hypothesis that $e(I)$ is free
of $d$-th powers implies that  $I$ is
projectively full.

Notice that for $(R,M)$ as in  Example \ref{5.2}, we have $e(M) =
k^2$.  In particular, for each integer $k \ge 2$,  Example \ref{5.2}
establishes the existence of a normal local domain $(R,M)$ of
altitude two such that $M$ is not projectively full and $e(M) =
k^2$.

}
\end{rema}

\section{Rational singularities and the closed fiber.}

  Lipman in his paper
\cite{Lipman.Rational.Sing} extends  Zariski's theory of complete
ideals in a  regular local domain of altitude two  to the case where
$R$ is a  normal
local domain of altitude two that has a rational singularity. Lipman proves that
$R$ has unique factorization of complete ideals if and only if the
completion of $R$ is a UFD. For $R$ having this property, it follows
that $\mathbf P(I)$ is projectively full for each nonzero proper
ideal $I$, e.g. this is the case if  $R = F[[x, y, z]]$,  where
$z^2 + y^3 + x^5 = 0$ and $F$ is a field of characteristic zero.

Let $(R,M)$ be a normal local domain of altitude two. G\"ohner
proves that if $R$ has a rational singularity, then the set of
complete
 asymptotically irreducible ideals associated to a prime
 $R$-divisor $v$
 consists of the powers of an ideal $A_v$ which is uniquely
 determined
 by $v$. In our terminology, this says that if $I$ is a nonzero
 proper
 ideal of $R$ having only one Rees valuation ring, then
 $\mathbf P(I)$ is
 projectively full. G\"ohner's proof involves choosing a
 desingularization $f : X \to \Spec R$ such that $v$ is centered
 on a
 component $E_1$ of the closed fiber on $X$.

 Let $E_2, \ldots, E_n$ be the other
components of the closed fiber on $X$.
Let $E_X$ denote the group of divisors having the form
$\sum_{i = 1}^nn_iE_i$,
with $n_i \in \Z$, and consider
$$E_X^+ = \{D \in E_X \, | \, D \ne 0 \text{ and } (D \cdot
E_i) \le 0  \text{ for all } 1 \le i \le n \}, \mbox{ and }$$
$$E_X^{\#} = \{ D \in E_X \, |\, D \ne 0 \text{ and } O(-D) \text{
is generated  by its sections over } X \}.
$$

Lipman shows $E_X^{\#} \subseteq E_X^+$ and  equality holds  if
$R$ has a rational
singularity. Also, if $D = \sum_in_iE_i \in E_X^+$, then
negative-definiteness
of the intersection matrix $(E_i \cdot E_j)$ \,\, implies $n_i \ge
0$ for all $i$.

For if $D \in E_X^+$ and $D = A - B$, where $A$ and $B$ are
effective, then
$(A-B \cdot B) \le 0$ and $(A \cdot B) \ge 0$ imply $(B \cdot B)
\ge 0$, \,\, so $B = 0$.

Let $v = v_1, v_2, \ldots, v_n$ denote the discrete valuations
corresponding
to $E_1, \ldots, E_n$.  Associated with  $D = \sum_in_iE_i \in
E_X^{\#}$ one defines  the
complete $M$-primary ideal
$$
I_D = \{ r \in R \, |\, v_i(r) \ge n_i
\text{ for } 1 \le i \le n \}.
$$

This sets up  a one-to-one correspondence between elements of
$E_X^{\#}$
and complete $M$-primary ideals that generate invertible
$O_X$-ideals.

 Lipman suggested to us the following
proof that $\mathbf P(I)$ is projectively   full
for each complete $M$-primary ideal $I$ if $R$ has a rational
singularity.
Fix a desingularization
$$f : X \to \Spec R
$$ such that $I$ generates an invertible $O_X$-ideal and let $D =
\sum_in_iE_i \in E_X^{\#}$ be the divisor associated to $I$. Let $g
= \gcd \{ n_i \} $. Since $E^+ = E^{\#}$, it follows that $(1/g) D
\in E^{\#}$.

The ideals $J \in \mathbf P(I)$ correspond to divisors in $E^{\#}$
that are integral
multiples of $(1/g)D$. Thus if $K$ is the complete $M$-primary ideal
associated to $(1/g)D$, then each $J \in \mathbf P(I)$ is the
integral closure of
a power of $K$, so $\mathbf P(I)$ is projectively full.

\section{Integral extensions.}

Perhaps the main unresolved question related to projectively full
ideals is the following:

\begin{ques} \label{7.1} {\em
Let $I$ be a nonzero ideal of a Noetherian domain $R$.
 Does there always exist a finite integral extension domain $B$ of $R$
such that $\mathbf P(IB)$ is projectively full?

}
\end{ques}

An affirmative answer to Question \ref{7.1} is obtained in
\cite[Theorem 2.5]{CHRR3} if the ideal $I$ has a Rees integer that
is a unit of $R$. More precisely, let $I = (b_1, \ldots, b_g)R$ be a
nonzero regular ideal of the Noetherian ring $R$, let $R_g = R[X_1,
\ldots, X_g]$ and $K = (X^c_1 - b_1, \ldots, X^c_g - b_g)R_g$, where
$c$ is a positive integer. Then $A = R_g/K$ is a finite free
integral extension ring of rank $c^g$ of $R$. Let $x_i = X_i \mod K$.
Then $J = (x_1, \ldots, x_g)A$ is such that $(IA)_a = (J^c)_a$, so
$IA$ and $J$ are projectively equivalent.

\begin{theo}
Let $I = (b_1, \ldots, b_g)R$ be as above and let $(V,N)$ be a  Rees
valuation ring of $I$. Assume  that  $b_iV = IV = N^c$ for each $i$
with $1 \le i \le g$,  and that $c$ is a unit of $R$. Then the
finite free  integral  extension ring $A = R_g/K$ is such that $J = (x_1,
\ldots, x_g)A$ is projectively  full. Hence $\mathbf P(IA)$ is
projectively full. If $R$ is an integral domain and $z$ is a minimal
prime of $A$,  then $B = A/z$ is an integral extension domain of $R$
such that $\mathbf P(IB)$ is projectively full.

\end{theo}

Let $I$ be a regular proper ideal in a Noetherian ring $R$ and let
$e$ be the least common multiple of the Rees integers of $I$.
Consider the Rees ring $R[t^{-1}, It]$ and let $B$ denote the
integral closure of  $R[t^{-1}, It][t^{-\frac{1}{e}}]$ in
$R[t^{-\frac{1}{e}}, t^{\frac{1}{e}}]$.
 Shiroh Itoh in \cite[page
392]{It} proves that $t^{-\frac{1}{e}}B$ is a radical ideal, and
$B/t^{-\frac{1}{e}}B$ is a reduced graded ring.

This result of Itoh motivated us to ask about the existence of an
integral extension ring $A$ of $R$ in which there exists a radical ideal
that is projectively equivalent to $IA$. The following two result
are established in \cite{HRR}.

\begin{theo}
\label{Itoh} Let  $I$  be a regular proper ideal in a Noetherian
ring  $R$ and let  $b_1,\dots,b_g$  be regular elements in  $R$ that
generate  $I$. Let  $m$  be a positive integer and let
$$ K = ({X_1}^m-b_1,\dots,{X_g}^m-b_g)R[X_1,\dots,X_g].
$$ For  $i$ $=$
$1,\dots,g$  let $x_i$  be the residue class of  $X_i$  modulo  $K$,
let  $A$  $=$ $R[x_1,\dots,x_g]$, and let  $J$ $=$
$(x_1,\dots,x_g)A$.  Then:

\noindent {\bf{(\ref{Itoh}.1)}} $A$  is a finite free integral
extension ring of  $R$ and  $(J^m)_a$ $=$ $(IA)_a$, so $J$  is
projectively equivalent to  $IA$.

\noindent {\bf{(\ref{Itoh}.2)}} $(J^k)_a \cap R$ $=$
$I_{\frac{k}{m}}$ ($=$ $\{x \in R \mid \overline {v}_I(x) \ge \frac{k}{m}\}$)
for all positive integers  $k$.

\noindent {\bf{(\ref{Itoh}.3)}} If  $e$  is the least common
multiple of the Rees integers of  $I$  and  $m$ $=$ $e$, then  $J_a$
is a radical ideal.

\noindent {\bf{(\ref{Itoh}.4)}} If  $R$  is an integral domain and
if  $z$  is a minimal prime ideal in  $A$, then $(J+z)/z$ is
projectively equivalent to  $(IA+z)/z$, and if  $m$ $=$ $e$  as in
(\ref{Itoh}.3), then  $(J_a+z)/z$  is a radical ideal.
\end{theo}

\begin{theo} \label{7.4}
\label{radical} Let  $I$  be a regular proper
ideal in a Noetherian
ring  $R$, let  $V_1,\dots,V_n$ be the Rees
valuation rings of  $I$,
and let  $e_1,\dots,e_n$ be the
corresponding Rees integers of  $I$.
Assume that the least common multiple
$e$ of  $e_1, \ldots, e_n$  is
a unit in  $R$ and that there exists a regular element $b$  in  $I$
such that  $bV_i$ $=$ $IV_i$  for  $i$ $=$ $1,\dots,n$.

\noindent {\bf{(\ref{radical}.1)}} Let  $b_1,\dots,b_g$  be regular
elements in  $R$  that generate  $I$, let $A$ $=$ $R[x_1,\dots,x_g]$
($=$ $R[X_1,\dots,X_g]/({X_1}^e-b_1,\dots,{X_g}^e-b_g)$), and let
$J$ $=$ $({x_1},\dots, {x_g})A$. Then  $A$  is a finite free
integral extension ring of  $R$, $J_a$ $=$ $\Rad(IA)$  is a
projectively full radical ideal that is projectively equivalent to
$IA$, $(IA)_a$ $=$ $(J^e)_a$, and all Rees integers of  $J$  are
equal to one.

\noindent {\bf{(\ref{radical}.2)}} If  $R$  is a Noetherian domain,
then let  $b_1$ $=$ $b$, $b_2,\dots,b_g$ be a basis of  $I$, for $i$
$=$ $1,\dots,g$ fix an  $e$-th root  ${b_i}^{1/e}$ of  $b_i$ in an
algebraic closure of the quotient field of  $R$, and let $A$ $=$
$R[{b_1}^{1/e},\dots,{b_g}^{1/e}]$ and  $J$ $=$
$({b_1}^{1/e},\dots,{b_g}^{1/e})A$. Then  $A$  is a Noetherian
domain that is a finite integral extension of  $R$, $\Rad(IA)$ $=$
$J_a$  is a projectively full radical ideal that is projectively
equivalent to  $IA$, $(IA)_a$ $=$ $(J^e)_a$, and all Rees integers
of  $J$  are equal to one.
\end{theo}

We have the following corollary to Theorem \ref{7.4}

\begin{coro}
\label{great} Let  $I$  be a regular nonzero ideal in a Noetherian
ring  $R$  and assume that $R$  contains the field of rational
numbers.  Then there exists a finite free integral extension ring
$A$ of  $R$  that contains a projectively full radical ideal $J$
that is projectively equivalent to  $IA$  and whose Rees integers
are all equal to one. If  $R$  is an integral domain, then for each
minimal prime ideal  $z$  in  $A$ the ideal $(J+z)/z$  in $A/z$ is a
projectively full radical ideal that is projectively equivalent to
$(IA+z)/z$, and all Rees integers of $(J+z)/z$  are equal to one.
\end{coro}

The following example indicates difficulties in establishing a
result similar to Theorem \ref{7.4} without the hypothesis that $e$
is a unit of $R$.

\begin{exam}
\label{bad} {\em Let  $D$ $=$ $k[y]$, where  $k$  is a field of
characteristic two and  $y$  is an indeterminate, and let  $b$ $=$
$y^2(y-1)$. Then  $D$  is a PID, and  $I$ $=$ $bD$  has the two Rees
integers  $e_1$ $=$ $1$  (for its Rees valuation ring  $D_{p_1}$,
where  $p_1$ $=$ $(y-1)D$)  and  $e_2$ $=$ $2$  (for its Rees
valuation ring  $D_{p_2}$, where  $p_2$ $=$ $yD$), so the least
common multiple of  $e_1,e_2$  is  $e$ $=$ $2$. Let  $E$  be the
integral closure of $D[x]$  in its quotient field, where  $x$ $=$
$\sqrt b$. Then  $E$  is a Dedekind domain and the Rees integer of
$xE$  with respect to  $E_{yE}$  is 2, so  $xE$  is not a radical
ideal of $E$. }
\end{exam}

\begin{flushleft}
Department of Mathematics, North Dakota State University, Fargo,
North Dakota 58105-5075 {\em E-mail address:
catalin.ciuperca@ndsu.edu}

\vspace{.15in}

Department of Mathematics, Purdue University, West Lafayette,
Indiana 47909-1395 {\em E-mail address: heinzer@math.purdue.edu}

\vspace{.15in}

Department of Mathematics, University of California, Riverside,
California 92521-0135 {\em E-mail address: ratliff@math.ucr.edu}

\vspace{.15in}

Department of Mathematics, University of California, Riverside,
California 92521-0135 {\em E-mail address: rush@math.ucr.edu}

\end{flushleft}
\end{document}